\magnification=1200
\input amstex

\documentstyle{amsppt}
\NoBlackBoxes

\pagewidth{13cm} \pageheight{18cm}
\def\leftitem#1{\item{\hbox to\parindent{\enspace#1\hfill}}}
\def\Vir{\operatorname{Vir}}
\def\rank{\operatorname{rank}}
\def\supp{\operatorname{supp}}
\def\dim{\operatorname{dim}}
\def\dim{\operatorname{dim}}
\def\HVir{\operatorname{HVir}}

\def\a{\alpha} \def\d{\delta}
\def\b{\beta}

\def\N{\Bbb N}
\TagsOnRight

\def\C{\Bbb C}

\def\Z{\Bbb Z}
\def\Q{\Bbb Q}

\def\proclaim#1{\par\vskip.25cm\noindent{\bf#1. }\begingroup\it}
\def\endproclaim{\endgroup\par\vskip.25cm}

\rightheadtext{ Harish-Chandra Modules Over $\HVir[\Q]$}
\leftheadtext{X. Guo, X Liu and K. Zhao}

\topmatter
\title Harish-Chandra modules over the $\Q$ Heisenberg-Virasoro Algebra
\endtitle
\author    Xiangqian Guo, Xuewen Liu and Kaiming Zhao
\endauthor


\abstract In this paper, it is proved that all irreducible
Harish-Chandra modules over the $\Q$ Heisenberg-Virasoro algebra are
of intermediate series (all weight spaces are $1$-dimensional).
\endabstract

\endtopmatter

\document

\subhead 1. \ \ Introduction
\endsubhead
\medskip

The Heisenberg-Virasoro algebra is the universal central extension
of the Lie algebra of differential operators on a circle of order at
most one. It contains the classical Heisenberg algebra and the
Virasoro algebra as subalgebras. The structure of their irreducible
highest weight modules was studied in [ACKP, B1]. Irreducible
Harris-Chandra modules over the Heisenberg-Virasoro algebra were
totally classified in [LZ]. They are either highest weight modules,
lowest weight modules, or modules of intermediate series. The
representation theory of the Heisenberg-Virasoro algebra is closely
related to those of other Lie algebras, such as the Virasoro algebra
and toroidal Lie algebras, see [B2, FO, JJ]. For other results on
generalized Heisenberg-Virasoro algebras, please see [FO, SS, SJ]
and references therein.
\smallskip

Recently, some authors introduced generalized Heisenberg-Virasoro
algebras and started to study their representations (see [LJ, SS]).
In the present paper, we give the classification of irreducible
Harish-Chandra modules over the $\Q$ Heisenberg-Virasoro algebra.
Similar to the $\Q$ Virasoro algebra case [Ma], they are only
modules of intermediate series. The main idea in our proof (Lemma
3.1) is similar to those used in [Ma] and [GLZ1].
\smallskip

In this paper we always denote by $\Z, \Q, \C$ the sets of integers,
rational and complex numbers, respectively.  Now we   give the
definitions of generalized Heisenberg-Virasoro algebras and the $\Q$
Heisenberg-Virasoro algebra.
\medskip

{\bf Definition 1.1.} Suppose that $G$ is an additive subgroup of
$\C$. The generalized Heisenberg-Virasoro algebra $\HVir[G]$ is a
Lie algebra over $\C$ with a basis:
$$\{d_g,I(g),C_{D},C_{DI},C_I| g\in G\},$$ subject to the Lie brackets given
by: $$[d_g,d_h]=(h-g)d_{g+h}+\d_{g,-h}\frac{g^3-g}{12}C_D,\eqno
(1.1)$$
$$[d_g,I(h)]=hI(g+h)+\d_{g,-h}(g^2+g)C_{DI},\eqno (1.2)$$
$$[I(g),I(h)]=g\d_{g,-h}C_I,\eqno (1.3)$$
$$[\HVir[G],C_D]=[\HVir[G],C_{DI}]=[\HVir[G],C_I]=0.\eqno (1.4)$$

It is easy to see that $\HVir[G]\simeq \HVir[G']$ iff there exists
nonzero $a\in\C$ such that $G=aG'$.  When $G=\Z$, $\HVir[\Z]=\HVir$
is the classical Heisenberg-Virasoro algebra. When $G=\Q$, we call
$\HVir[\Q]$ the {\it $\Q$ Heisenberg-Virasoro algebra}, which is the
main object of this paper. An $\HVir[G]$ module $V$ is said to be
trivial if $\HVir[G]V=0$, and we denote the $1$-dimensional trivial
module by $T$.
\smallskip

The modules of intermediate series $V(\a,\b;F)$ over $\HVir[G]$ are
defined as follows for every $\a,\b,F\in\C$. As a vector space over
$\C$, $V(\a,\b;F)$ has a basis $\{v_g|g\in G\}$ and the actions are:
$$d_gv_h=(\a+h+g\b)v_{g+h}, \,\, I(g)v_h=Fv_{g+h},\eqno (1.5)$$
$$C_Dv_g=0,\,\,  C_Iv_g=0,\,\, C_{DI}v_g=0,\,\, \forall\, g,h\in G.\eqno (1.6)$$
\smallskip

It is well known that $V(\a,\b;F)\cong V(\a+g,\b;F)$ for any
$\a,\b,F\in\C$ and $g\in G$. So we always assume that $\a=0$ when
$\a\in G$.    It is also easy to see that $V(\a,\b;F)$ is reducible
if and only if $F=0, \a=0$, and $\b\in\{0,1\}$. The module
$V(0,0;0)$ has a $1$-dimensional submodule $T$ and $V(0,0;0)/T$ is
irreducible; $V(0,1;0)$ has codimension $1$  irreducible submodule.
We denote the unique nontrivial irreducible sub-quotient of
$V(\a,\b;F)$ by $V'(\a,\b;F)$. It is easy to see that
$V'(0,0;0)\cong V'(0,1;0)$.

{\bf Remark.} It is not hard to verify that, for $F\ne0$,
$V(\a,\b;F)\cong V(\a',\b';F')$ iff $\a-\a'\in G$ and $(\b, F)=(\b',
F')$.
\smallskip

Our main theorem is the following:
\smallskip

\proclaim{Theorem 1.2} Suppose that $V$ is an irreducible nontrivial
Harish-Chandra module over $\HVir[\Q]$. Then $V$ is isomorphic to
 $V'(\a,\b;F)$ for suitable $\a,\b, F\in\C$.\endproclaim

The paper is organized as follows. In Section 2, we collect some
known results on the Virasoro algebra and on the Heisenberg-Virasoro
algebra for later use. In section 3, we give the proof of the main
theorem.

%
%
%

 \subhead 2. \ \ Preliminaries
\endsubhead
\smallskip

Since generalized Heisenberg-Virasoro algebras are closely related
to generalized Virasoro algebras, we first recall some results on
generalized Virasoro algebras.
\smallskip

{\bf Definition 2.1} Let $G$ be a nonzero additive subgroup of $\C$.
The generalized Virasoro algebra $\Vir[G]$ is a Lie algebra over
$\C$ with a basis: $\{d_g,C_D| g\in G\}$ and subject to the Lie
brackets given by:
$$[d_g,d_h]=(h-g)d_{g+h}+\d_{g,-h}\frac{g^3-g}{12}C_D,\,\,\,\,\,\,\, [\Vir[G],C_D]=0.\eqno
(2.1)$$
\smallskip

Notice that $\Vir[G]$ is a subalgebra of $\HVir[G]$. When $G=\Z$,
$\Vir[\Z]=\Vir$ is the classical Virasoro algebra and when $G=\Q$,
$\Vir[\Q]$ is called the $\Q$ Virasoro algebra.
\smallskip

A module of intermediate series $V(\a,\b)$ has a $\C$-basis
$\{v_g|g\in G\}$ and the $\Vir$-actions are:
$$d_gv_h=(\a+h+g\b)v_{g+h},\,\,\,\,\,\,\,\,\,\,C_Dv_g=0,\,\,\,\,\,\,\,\,\,\, \forall g,h\in G.\eqno (2.2)$$
\smallskip

It is well known that $V(\a,\b)$ is reducible if and only if $\a\in
G$ and $\b\in\{0,1\}$. The unique nontrivial irreducible
sub-quotient of $V(\a,\b)$ is denoted by $V'(\a,\b)$. It is also
well known that $V(\a,\b)\cong V(\a+g,\b)$ for any $\a,\b\in\C$ and
any $g\in G$.  So we always assume that $\a=0$ if $\a\in G$.
\smallskip

Mazorchuk [Ma] classified  irreducible Harish-Chandra modules over
$Vir[\Q]$:

\proclaim{Theorem 2.2} Any irreducible Harish-Chandra module over
$\Vir[\Q]$ is isomorphic to $V'(\a,\b)$ for suitable
$\a,\b\in\C$.\endproclaim
\medskip

Now we consider the $\Q$ Heisenberg-Virasoro algebra $\HVir[\Q]$,
which contains a classical Heisenberg-Virasoro algebra
$\HVir[\Z]=\HVir$. The classification of irreducible Harish-Chandra
modules over $\HVir$ was given in [LZ]:

\proclaim{Theorem 2.3} Any irreducible Harish-Chandra module over
$\HVir$ is isomorphic to either a highest weight module, a lowest
weight module, or   $V'(\a,\b;F)$  for suitable $\a,\b,
F\in\C$.\endproclaim
\medskip

 \subhead 3. \ \ Proof of Theorem 1.2
\endsubhead
\smallskip

We decompose $\Q$ as the union of a series of rings
$\Q_k=\{\frac{n}{k!}\,|\, n\in\Z\},\,k\in\N$. Then we can viewe each
$\HVir[\Q_k]$ as a subalgebra of $\HVir[\Q]$ naturally, and thus
$\HVir[\Q]=\bigcup_{k\in\N}\HVir[\Q_k]$. Clearly, $\HVir[\Q_k]\simeq
\HVir$ for any $k\in\N$. For convenience, we write
$U(\Q)=U(\HVir[\Q])$ and $U(\Q_k)=U(\HVir[\Q_k])$ for short. Denote
$U(\Q)_q=\{u\in U(\Q)| [d_0,u]=qu\}$ and $U(\Q_k)_q=\{u\in U(\Q_k)|
[d_0,u]=qu\}$, for any $q\in\Q$.

\proclaim{Lemma 3.1}Suppose that $N$ is a finite dimensional
irreducible $U(\Q)_{0}$-module, then there is some $k\in\N$ such
that $N$ is an irreducible $U(\Q_{k})_{0}$-module.
\smallskip

\demo{Proof} There is an associative algebra homomorphism $\Phi$:
$U(\Q)_{0}\longrightarrow gl(N)$, where $gl(N)$ is the general
linear associative algebra of $N$.
\smallskip

Since $N$ and hence $gl(N)$ are both finite dimensional, then
$U(\Q)_{0}/\ker(\Phi)$ is finite dimensional. Take $y_1$,
$y_{2}$,$\ldots$, $y_{m}\in U(\Q)_{0}$ such that $\overline{y_{1}}$,
$\overline{y_{2}}$, $\cdots$, $\overline{y_{m}}$ become a basis of
$U(\Q)_{0}/\ker(\Phi)$. Then there is some $k\in \N$ such that
$y_{1}$, $y_{2}$,$\ldots$, $y_{m}$ are all in $U(\Q_{k})_{0}$.
\smallskip

We claim that $N$ is an irreducible $U(\Q_{k})_{0}$-module. Let $M$
be a proper $U(\Q_{k})_{0}$-submodule of $N$. For any $y\in
U(\Q)_{0}$, there is some $y_{0}\in \ker(\Phi)$ and  $a_{i}\in \C$,
such that $y=y_{0}+\sum_{i=1}^m a_{i}y_{i}$. Thus
$yM=(y_{0}+\sum_{i=1}^m a_{i}y_{i})M\subset(\sum a_{i}y_{i})M\subset
M$, that is, $M$ is a $U(\Q)_{0}$-submodule of $N$, forcing $M=0$.
Thus $N$ is irreducible over $U(\Q_{k})_{0}$.\qed
\enddemo
\smallskip

Now we fixe a nontrivial irreducible Harish-Chandra module $V$ over
  $\HVir[\Q]$. Then, there exists
some $\a\in\C$ such that $V=\bigoplus_{q\in \Q}V_q$, where
$V_q=\{v\in V\,|\,d_0v=(\a+q)v\}$. We define the support of $V$ as
$\supp V=\{q\in\Q| V_q\neq 0\}$.
%
%
%

\proclaim{Lemma 3.2}$\supp V=\Q$ or $\Q\setminus\{-\a\}$ and $\dim
V_{q}=1$, $\forall \, q\in \supp V$.
\smallskip

\demo{Proof} First view $V$ as a $Vir[\Q]$ module, by Theorem 2.2,
then $\dim V_p=\dim V_q, \,\forall\, p, q\in\Q\setminus\{-\a\}$.
Particularly,  $V$ is a uniformly bounded module, i.e., the
dimensions of all weight spaces are bounded by a positive integer.
\smallskip

Suppose that $\dim V_{q}>1$ for some $q\in \supp V$. Then it is easy
to see that $V_q$ is an irreducible $U[\Q]_{0}$-module. By Lemma
3.1, there is some $k\in \N$ such that $V_{q}$ is an irreducible
$U(Q_{k})_{0}$-module.
\smallskip

We consider the $\HVir[\Q_k]$-module $W=U(\Q_{k})V_{q}$, which is
uniformly bounded. Then by Theorem 2.3, there is a composition
series of $\HVir[\Q_{k}]$-modules:
$$0=W^{(0)}\subset W^{(1)}\subset W^{(2)}\subset\ldots
W^{(m)}=W.$$ Each factor $V^{(i)}/V^{(i-1)}$ is either trivial or of
intermediate series, and in both cases all nonzero weight spaces are
$1$-dimensional.

\smallskip

Take the first $V^{(i)}$ such that $\dim V^{(i)}_q\neq0$. Then $\dim
V^{(i)}_{q}=\dim(V^{(i)}/V^{(i-1)})_{q}=1$. But on the other hand,
$V^{(i)}_{q}$ is a $U(\Q_{k})_{0}$-module, i.e., a nontrivial proper
$U(\Q_{k})_{0}$-submodule of $V_{q}$, contradiction. Thus $\dim
V_{q}=1$, $\forall \, q\in \supp V$.

Since $V$ is nontrivial, we must have some $q\in\supp
V\setminus\{-\a\}$. Our result follows since $\dim V_p=\dim
V_q=1,\,\,\forall\, p\in\Q\setminus\{-\a\}$.\qed
\enddemo
\medskip

Now we can give the proof of our main theorem:
\smallskip
%
%

\demo{Proof of Theorem 1.2} By Lemma 3.2, we know that $\supp V=\Q$
for some $\a\in\C$ or $\supp V=\Q\setminus\{-\a\}$. Denote
$V^{(0)}=0$ and $V^{(k)}=\bigoplus_{q\in \Q_k} V_{q}$ for all $k\geq
1$. Then $\supp V^{(k)}=\{q\in\Q| V^{(k)}_q\neq 0\}=\supp V\cap\Q_k$
for all $k\geq 1$. We  have a vector space filtration of $V$:
$$0=V^{(0)}\subset V^{(1)}\subset V^{(2)}\subset..\subset V^{(k)}\subset...\subset V \,\,\text{and}\,\,V=\bigcup_{k=0}^{\infty}V^{(k)}.$$
\smallskip

It is clear that each $V^{(k)}$ can be viewed as an
$\HVir[\Q_k]$-module, and each $\HVir[\Q_k]$ is isomorphic to
$\HVir$.
\medskip

If $-\a\notin\supp V$, then each $V^{(k)}$ is irreducible over
$\HVir[\Q_k]$ for any $k\in\Z$ by Lemma 3.2 and Theorem 2.3.
\smallskip

If $-\a\in\supp V$, we have assumed that $\a=0$ in this case.  Then
we have $u_p\in U(\HVir[\Q])_p$ and $u_q\in U(\HVir[\Q])_q$ for some
$0\neq p,q\in\Q$ such that $u_pV_0=V_p$ and $u_qV_{-q}=V_0$, since
$V$ is an irreducible $\HVir[\Q]$-module. Then there is some $k_0$
such that $u_p,u_q\in U(\HVir[\Q_{k_0}])$, thus $V^{(k_0)}$ is
irreducible as an $\HVir[\Q_{k_0}]$-module.
\medskip

Since $V^{(k_0)}$ is irreducible over $\HVir[\Q_{k_0}]$, then
$V^{(k)}$ is irreducible over $\HVir[\Q_{k}]$ for any $k\geq k_0$.
By Theorem 2.3, we see that $C_DV=C_IV=C_{DI}V=0$, that $I(0)$ acts
as a scalar $F\in\C$, and that there are some $\a_k, \b_k\in\C$ such
that $V^{(k)}\cong V'(\a_k,\b_k;F)$ as modules over $\HVir[\Q_{k}]$
with $k\ge k_0$. Write $\a=\a_{k_0}$ and $\b=\b_{k_0}$ for short.
\medskip

If $F=0$, $V$ is an irreducible $\Vir(\Q)$-module and $I(p)V=0$ for
any $p\in\Q$. Theorem 2.1 follows from Theorem 2.2. Next we assume
that $F\ne0$.
\smallskip

Now we need to prove that $V\cong V'(\a,\b;F)$ as modules over
$\HVir[\Q]$, i.e., we can  choose a basis of $V$ such that (1.5) and
(1.6) hold. We proceed by choosing a basis of each $V^{(k)}$
inductively, such that (1.5) and (1.6) hold when one only consider
the actions of $\HVir[\Q_k]$, and that the basis of each $V^{(k+1)}$
is the extension of the basis of $V^{(k)}$ for any $k\geq k_0$.
\smallskip

Naturally, we have such a basis for $V^{(k_0)}$. Now suppose that we
have such bases for $V^{(k)},\,\forall\, k<m$ for some $m>k_0$. That
is, we have a basis $\{v_q\in V_q\,|\, q\in \supp V^{(m-1)}\}$ such
that $d_pv_q=(\a+q+p\b)v_{p+q}$ and $I(p)v_q=Fv_{p+q}, p,q\in \supp
V^{(m-1)}$. Now we consider $V^{(m)}$.
\medskip

Let $\Phi$ be the canonical isomorphism $\HVir\,\,\longrightarrow
\HVir[\Q_{m}],$ defined by $\Phi(d_n)=m!d_{\frac{n}{m!}}$ and
$\Phi(I(n))=m!I(\frac{n}{m!}),\,\ \forall\, n\in\Z$.
\smallskip

Then $V^{(m)}$ can be viewed as an $\HVir$ module via $\Phi$ and is
isomorphic to $V(\a',\b';F')$ as modules over $\HVir$ for suitable
$\a',\b',F'\in\C$, by Theorem 2.3. Then there is a basis
$\{v'_q|q\in \supp V^{(m)}\}$ of $V^{(m)}$ satisfying:
$$(m!d_{\frac{i}{m!}})v'_{\frac{k}{m!}}=(\a'+k+i\b')v'_{\frac{k+i}{m!}}, \,\,
(m!I(\frac{i}{m!}))v'_{\frac{k}{m!}}=F'v'_{\frac{k+i}{m!}},\,\,
\forall \frac{i}{m!},\frac{k}{m!}\in\supp V^{(m)}.$$
\smallskip

That is:
$$d_{\frac{i}{m!}}v'_{\frac{k}{m!}}=(\frac{\a'}{m!}+\frac{k}{m!}+\frac{i}{m!}\b')v'_{\frac{k+i}{m!}}, \,\,
I(\frac{i}{m!})v'_{\frac{k}{m!}}=\frac{F'}{m!}v'_{\frac{k+i}{m!}},$$
or that,
$$d_{p}v'_{q}=(\frac{\a'}{m!}+q+p\b')v'_{p+q}, \,\,
I(p)v'_{q}=\frac{F'}{m!}v'_{p+q},\,\forall\, p,q\in\supp V^{(m)}.$$
\medskip

Taking $p=0$ and $q\in\supp V^{(m-1)}\subset\supp V^{(m)}$, we then
have that $d_{0}v'_{q}=(\frac{\a'}{m!}+q)v'_{q}$ and
$I(0)v'_{q}=\frac{F'}{m!}v'_{q}$, which indicates that $\a'=m!\a$
and $F'=m!F$.
\smallskip

Assume that $v'_q=c_qv_q,\,\forall\, q\in\supp V^{(m-1)}$. Compare
$I(p)v'_{q}=Fv'_{p+q}$ and $I(p)v_{q}=Fv_{p+q}$. We see that $c_p$
is independent to $p$. Then we may choose $v'_p=v_p$ for all
$p\in\supp V^{(m-1)}$. By the remark in Sect.2 we see that $\b'=\b$.
Denote $v_q=v'_q$ for all $q\in\supp V^{(m)}\setminus\supp
V^{(m-1)}$. Thus the basis $\{v_q|q\in\supp V^{(m)}\}$ is an
extension of $\{v_p|p\i\supp V^{(m-1)}\}$ such that
$d_pv_q=(\a+q+p\b)v_{p+q}$ and $I(p)v_q=Fv_{p+q},\,\forall\,
p,q\in\supp V^{(m)}.$
\medskip

By induction, we can have a basis $\{v_q\,|\,q\in\supp V\}$ for $V$
such that $d_pv_q=(\a+q+p\b)v_{p+q}$, $I(p)v_q=Fv_{p+q},\,\forall\,
p,q\in\supp V$ and $C_DV=C_IV=C_{DI}V=0$. That is, $V\cong
V'(\a,\b;F)$. This completes the proof. \qed
\enddemo
\smallskip

Recall from [K] that the {\bf rank} of an additive subgroup $G$ of
$\C$, denoted by {\bf rank$(G)$}, is the maximal number $r$ with
$g_1,\cdots,g_r\in G\setminus\{0\}$ such that $\Z g_1+...+\Z g_r$ is
a direct sum. If such an $r$ does not exist, we define
$\rank(G)=\infty$.
\smallskip

Now we can see that the proof of Theorem 3.3 is also valid for $G$
being an infinitely generated additive subgroup of $\C$ of rank $1$.
Thus, we have the following:

\smallskip
\proclaim{Theorem 3.3}Let $G$ be an infinitely generated additive
subgroup of $\C$ with rank $1$. Then any nontrivial irreducible
Harish-Chandra module over $\HVir[G]$ must be isomorphic to
$V'(\a,\b;F)$ for suitable $\a,\b, F\in\C$.\qed
\endproclaim

\vskip.3cm \Refs\nofrills{\bf REFERENCES}
\bigskip
\parindent=0.45in

\leftitem{[ACKP]} E. Arbarello, C. De Concini, V. G. Kac and C.
Procesi, Moduli spaces of curves and representation theory, {\it
Comm. Math. Phys.}, (1)117, 1-36(1988).

\leftitem{[B1]} Y. Billig,  Representations of the twisted
Heisenberg-Virasoro algebra at level zero, {\it Canad. Math. Bull.},
46(2003), no.4, 529-537.

\leftitem{[B2]} Y. Billig, category of modules for the full toroidal
Lie algebra, {\it Int. Math. Res. Not.}, 2006, Art. ID 68395, 46pp.

\leftitem{[BGLZ]}  P. Batra, X. Guo, R. Lu  and  K. Zhao, Highest
weight modules over the pre-exp-polynomial algebras, {\it J.
Algebra}, Vol.322, 4163-4180(2009).

\leftitem{[FO]} M. Fbbri and F. Okoh, Representation of Virasoro
Heisenbergs and Virasoro toroidal algebras, Canad. J. Math.,
51(1999), no.3, 523-545.

\leftitem{[GLZ1]} X. Guo, R. Lu and K. Zhao, Classification of
irreducible Harish-Chandra modules over the generalized Virasoro
algebra, Preprint.

\leftitem{[GLZ2]} X. Guo, R. Lu and K. Zhao, Simple Harish-Chandra
modules, intermediate series modules and Verma modules over the
loop-Virasoro algebra, Forum Math., in press.

\leftitem{[JJ]} Q. Jiang and C. Jiang, Representations of the
twisted Heisenberg-Virasoro algebra and the full toroidal Lie
algebras, to appear in {\it Algebra Colloq.}

\leftitem{[K]} I. Kaplansky, {\it Infinite abelian groups}, Revised
edition, The University of Michigan Press, Ann Arbor, Mich. 1969.

\leftitem{[KR]} V. G. Kac and K. A. Raina, ``Bombay lectures on
highest weight representations of infinite dimensional Lie
algebras,'' World Sci., Singapore, 1987.

\leftitem{[LJ]} D. Liu and C. Jiang, The generalized Heisenberg
Virasoro algebra, Preprint.

\leftitem{[LZ]} R. Lu and K. Zhao, Classification of irreducible
weight modules over the twisted Heisenberg-Virasoro algebra, to
appear in {\it Commun. Contemp. Math.}

\leftitem{[M]} O. Mathieu, Classification of Harish-Chandra modules
over the Virasoro algebra, {\it Invent. Math.} {\bf 107}(1992),
225-234.

\leftitem{[Ma]} V. Mazorchuk, Classification of simple
Harish-Chandra modules over $Q$ Virasoro algebra, {\it Math. Nachr.}
{\bf 209}(2000), 171-177.

\leftitem{[SJ]} R. Shen, C. Jiang, Derivation algebra and
automorphism group of the twisted Heisenberg-Virasoro algebra, to
appear in {\it Cmmm. Algebra}.


\leftitem{[SS]} R. Shen and Y. Su, Verma modules over the
generalized Heisenberg-Virasoro algebra, Comm. Algebra 36 (2008),
no. 4, 1464--1473.

%
%


\endRefs
\bigskip

\noindent Xiangqian Guo and Xuewen Liu: Department of Mathematics,
Zhengzhou university, Zhengzhou 450001, Henan, P. R. China.

\noindent Email: guoxq\@zzu.edu.cn and liuxw\@zzu.edu.cn
\medskip

\noindent Kaiming Zhao: Department of Mathematics, Wilfrid Laurier
University, Waterloo, ON, Canada N2L 3C5, and Academy of Mathematics
and System Sciences, Chinese Academy of Sciences, Beijing 100190, P.
R. China.

\noindent Email: kzhao\@wlu.ca

\vfill

\enddocument